# Asymmetric non-Gaussian effects in a tumor growth model with immunization [*]

Mengli Hao[1], Jinqiao Duan[2,3], Renming Song[4] and Wei Xu[1]
1. Department of Applied Mathematics
Northwestern Polytechnical University, Xi'an, 710129, China
E-mail: haomengli@gmail.com
2. Institute for Pure and Applied Mathematics
University of California Los Angeles, CA 90095, USA
E-mail: jduan@ipam.ucla.edu
3. Department of Applied Mathematics
Illinois Institute of Technology, Chicago, IL 60616, USA
E-mail: duan@iit.edu
4. Department of Mathematics
University of Illinois at Urbana-Champaign Urbana, IL 61801, USA
E-mail: rsong@math.uiuc.edu

### Abstract

The dynamical evolution of a tumor growth model, under immune surveillance and subject to asymmetric non-Gaussian $\alpha$-stable Lévy noise, is explored. The lifetime of a tumor staying in the range between the tumor-free state and the stable tumor state, and the likelihood of noise-inducing tumor extinction, are characterized by the mean exit time (also called mean residence time) and the escape probability, respectively. For various initial densities of tumor cells, the mean exit time and the escape probability are computed with different noise parameters. It is observed that unlike the Gaussian noise or symmetric non-Gaussian noise, the asymmetric non-Gaussian noise plays a constructive role in the tumor evolution in this simple model. By adjusting the noise parameters, the mean exit time can be shortened and the escape probability can be increased, simultaneously. This suggests that a tumor may be mitigated with higher probability in a shorter time, under certain external environmental stimuli.

[*]This work was done while Mengli Hao was visiting the Institute for Pure and Applied Mathematics (IPAM), Los Angeles, CA 90095, USA. This work was supported by the NSFC Grants 11172233, 10932009, 10971225 and 11028102, NSF Grant 1025422 and Simons Foundation grant 208236.

# 1 Introduction

In the recent years, more and more facts have illustrated the important influences of noise on dynamical systems. It is often assumed that the external noise is Gaussian. This arises due to the assumption that the external perturbation is the result of a large number of independent interactions of bounded strength [4, 8, 10, 14, 20]. However, this assumption is not always suitable to adequately interpret real data. For instance, when the fluctuations are abrupt pulses or extreme events, the Gaussian assumption is obviously not proper. In this case, it may be more appropriate to model the fluctuations by a process with heavy tails and discontinuous sample paths. A class of this kind of processes is the asymmetric $\alpha$-stable Lévy motion. Noises following symmetric and asymmetric $\alpha$-stable laws are abundant in nature and have been observed in various fields of science [3, 8, 17].

Moreover, due to the increasing number of people with tumor, cancer research has become a major challenge in medicine and biology. Because surgeries, chemotherapies and radiotherapies could bring great pain to patients and adversely affect patients' life, considerable attention has been paid to understanding immunotherapy, aiming to strengthen the body's own natural ability to combat cancer by enhancing the effectiveness of the immune system. One of the deterministic tested models representing the interactions between tumor tissue and immune system is via a "predator-prey" mechanism, in which tumor cells play the role of "preys" whereas the immune cells act as as "predators" [4, 12, 14, 16, 29].

A lot of studies have been conducted on the dynamical behaviors of this model driven by different noises, such as Gaussian or Brownian noise, colored noise, and fractional Gaussian noise [9, 10, 19, 20, 22, 27, 28]. These researchers have shown that environmental noises have great impact on the growth and extinction of tumor cells. Furthermore, the dynamical evolution behaves differently under different noises. Thus, it is important to gain deeper insight into the effects of various noises on this tumor evolution system.

In this paper, we focus on effects of the asymmetry of the non-Gaussian $\alpha$-stable Lévy noise on the dynamical behaviors of a tumor cell growth model with immunization. We evaluate the evolution of the tumor cell density, including the lifetime in the range between the tumor-free state and the stable tumor state, as well as the likelihood that the tumor cells become extinct, by numerically examining the mean exit time and the escape probability.

This paper is organized as follows. In Section 2, we introduce the asymmetric $\alpha$-stable Lévy processes. Section 3 describes the tumor growth model with immunization, under the influence of asymmetric $\alpha$-stable Lévy noise. In Section 4, we recall the integro-differential equations satisfied by the mean exit time $u$ and the escape probability $p$, and give numerical algorithms for solving these equations. We present numerical results in Section 5, then finish with conclusions in the final section.

## 2    Asymmetric $\alpha$−stable Lévy motions

A scalar Lévy motion is characterized by a drift parameter $\mu$, a non-negative diffusion constant $d$, and a non-negative Borel measure $\nu$, defined on $(\mathbb{R}^1, \mathcal{B}(\mathbb{R}^1))$ and concentrated on $\mathbb{R}^1 \setminus \{0\}$. The measure $\nu$ is called the Lévy jump measure and it has the following property:

$$\int_{\mathbb{R}^1 \setminus \{0\}} (y^2 \wedge 1)\, \nu(\mathrm{d}y) < \infty, \tag{1}$$

where $a \wedge b = \min\{a, b\}$. We call $(\mu, d, \nu)$ the generating triplet of the Lévy motion $L_t$.

The generator $A$ of the Lévy motion $L_t$ is defined by $A\varphi = \lim_{t \downarrow 0} \frac{\mathbb{P}_t \varphi - \varphi}{t}$, where $\mathbb{P}_t \varphi(x) = \mathbb{E}\varphi(L_t)$ and $\varphi$ is any function belonging to the domain of the operator $A$. Recall that the space $C_b^2(\mathbb{R}^1)$, consisting of $C^2$ functions with bounded derivatives up to order 2, is contained in the domain of $A$, and thus for every $\varphi \in C_b^2(\mathbb{R}^1)$ (see [1, 2])

$$A\varphi(x) = \mu \varphi'(x) + \frac{d}{2}\varphi''(x) + \int_{\mathbb{R}^1 \setminus \{0\}} [\varphi(x+y) - \varphi(x) - \mathbf{1}_{\{|y|<1\}}\, y\varphi'(x)]\, \nu(\mathrm{d}y). \tag{2}$$

In this paper, we consider a special scalar Lévy motion $L_t$ with the generating triplet $(0, d, \nu_{\alpha,\beta})$, for the diffusion coeffient $d \geq 0$, stability index $\alpha \in (0, 2)$, skewness parameter $\beta \in [-1, 1]$ and the Lévy jump measure $\nu_{\alpha,\beta}$.

More specifically, the jump measure is defined by [13, 23]:

$$\nu_{\alpha,\beta}(\mathrm{d}y) = \begin{cases} \dfrac{C_1 dy}{|y|^{\alpha+1}}, & y > 0, \\ \dfrac{C_2 dy}{|y|^{\alpha+1}}, & y < 0, \end{cases} \qquad (3)$$

with $\beta = (C_1 - C_2)/(C_1 + C_2)$, $C_1 = C_\alpha \frac{1+\beta}{2}$ and $C_2 = C_\alpha \frac{1-\beta}{2}$, where

$$C_\alpha = \begin{cases} \dfrac{\alpha(1-\alpha)}{\Gamma(2-\alpha)\cos(\frac{\pi\alpha}{2})}, & \alpha \neq 1, \\ \dfrac{2}{\pi}, & \alpha = 1. \end{cases}$$

This $L_t$ is a non-Gaussian process, although it has a Gaussian diffusion component described by the diffusion constant $d$. When $d = 0$, this is the well-known asymmetric $\alpha$-stable Lévy motion, and if additionally $\beta = 0$, it is the symmetric $\alpha$-stable Lévy motion.

Note that for a Lévy motion $L_t$, its characteristic function is [2]

$$\mathbb{E}[e^{i\lambda L_t}] = e^{-t\Psi(\lambda)}, \qquad t \geq 0, \lambda \in \mathbb{R}^1,$$

with $\Psi(\lambda)$ a function that does not depend on time $t$. In other words, $e^{-\Psi(\lambda)}$ is the characteristic function for $L_1$. The characteristic function of $L_t$ (at time $t$) is just the $t$−th power of the characteristic function of $L_1$ (at time 1).

For the $\alpha$−stable Lévy motion considered here ($d = 0$), $L_1$ is a $\alpha$−stable random variable with distribution $S_\alpha(1, \beta, 0)$; see [15, Ch. 1]. Figure 1 shows the probability density function (PDF) for the $\alpha$−stable random variable $L_1$, with various stability indices $\alpha$ and skewness parameters $\beta$.

# 3 Tumor growth model subject to an asymmetric Lévy noise

In this section, we consider a tumor growth model under immune surveillance. The reaction between the tumor tissues and the immune cells is based on a reaction scheme representative of the catalytic Michaelis-Menten scenario [10]. It can be explained as follows: Firstly, the tumor cells denoted by $X$ proliferate in two ways: one is the transformation of normal cells into neoplastic ones $X$ at a rate $\kappa$; the other is the replication of the tumor cells at a rate $\iota$. Then, the active cytotoxic cells (i.e. immune cells) $Y$ bind the tumor cells to the complex $Z$ with the kinetic constant $k_1$. Lastly, the complex $Z$ dissociates into immune cells and the dead or non-replicating tumor $P$ at a rate proportional to $k_2$. Schematically, the mechanisms above can be represented as follows:

$$\text{Normal Cells} \xrightarrow{\kappa} X, \tag{4a}$$

$$X \xrightarrow{\iota} 2X, \tag{4b}$$

$$X + Y \xrightarrow{k_1} Z \xrightarrow{k_2} Y + P. \tag{4c}$$

In order to construct a mathematical model, we can make the following assumptions based on biological principles. Firstly, because the transformation of normal cells into neoplastic ones originates from environmental carcinogenic agents rather than spontaneous endogenous somatic mutations, the average rate of this process is very low, compared with the rate of neoplastic cell replication [14,18] (Typical experimental values are: $\kappa$ the order of $10^{-17}$-$10^{-18}$ transformed cell/normal cell×day, $\iota$ =0.2-1.5 day$^{-1}$, $k_1$ =0.1-18 day$^{-1}$, $k_2$ =0.2-18 day$^{-1}$ as in [9,18]). Thus, we can ignore step (4a). Secondly, in the reaction process, $Y$ behaves like the enzymes in the Michaelis-Mentecn reaction, so we consider a conserved mass of enzymes $Y + Z = E = $ const. Besides, in the limit case, the production of X-type cells inhibited by a hyperbolic activation is the slowest process. Therefore, by the assumptions and the quasi-steady-state approximation, the kinetics can be simplified to an equivalent single variable differential equation [9,18]:

$$\frac{dx}{dt} = x(1 - \theta x) - \gamma \frac{x}{x+1}, \tag{5}$$

with the potential function

$$U(x) = -\frac{x^2}{2} + \frac{\theta x^3}{3} + \gamma x - \gamma \ln(x+1), \tag{6}$$

where $x$ is the normalized molecular density of tumor cells with respect to the maximum tissue capacity. And we use the following scaling formulas in the process of nondimensionalization:

$$x = \frac{k_1}{k_2} X, \theta = \frac{k_2}{k_1}, \gamma = \frac{k_1 E}{\iota}, t = \iota t'.$$

Taking into account the biological significance and convenience of discussion, we choose the parameter ranges as follows: $\theta < 1$, $0 < \gamma < \frac{(1+\theta)^2}{4\theta}$. In this case, the deterministic dynamical system Eq. (5) has two stable states and one unstable state (see [18, 20]). Namely, the potential function $U(x)$ has two minima: $x_1$, $x_3$ and one maximum $x_2$ (see Figure 2). That is to say, this system has three meaningful steady states:

$$x_1 = 0,$$
$$x_2 = \frac{1 - \theta - \sqrt{(1+\theta)^2 - 4\gamma\theta}}{2\theta}, \qquad (7)$$
$$x_3 = \frac{1 - \theta + \sqrt{(1+\theta)^2 - 4\gamma\theta}}{2\theta}.$$

Without random fluctuations, for a given initial condition, the system state will eventually converge to one of the two stable states: (i) the stable state $x_1 = 0$, called the tumor-free state (or the state of tumor extinction), where no tumor cells exist, (ii) the stable state $x_3$, called the stable tumor state, where the tumor cell density does not increase but keeps at a certain constant level.

However, from a biological point of view, the growth rate of tumor tissue is inevitably influenced by many environmental factors, such as temperature, radiations, chemical agents, the degree of vascularization of tissue, the supply of nutrients, the immunological state of the host and so on [9, 18]. Due to the limitations of the Gaussian noise, in this paper, we consider more general noise, an asymmetric Lévy noise to represent the environmental fluctuations. Our model is written as follows:

$$dX_t = f(X_t)dt + dL_t, \quad X(0) = x, \qquad (8)$$

where

$$f(X_t) = X_t(1 - \theta X_t) - \gamma \frac{X_t}{X_t + 1},$$

and $L_t$ is a Lévy process with the generating triplet $(0, d, \nu_{\alpha,\beta})$, i.e., drift coefficient 0, diffusion coefficient $d \geq 0$ and the Lévy jump measure $\nu_{\alpha,\beta}$.

Under the effects of environmental fluctuations, the number of tumor cells may fluctuate in the range between the tumor-free state and the stable tumor state denoted by $D = (x_1, x_3)$. In this paper, we concentrate on how the tumor density evolves in the range $D$ and discuss the impacts of the asymmetric Lévy noise on the time that the density of tumor cells stays in the range $D$ and the probability it exits $D$ from the left side, i.e., becoming tumor-free.

# 4 Mean exit time, escape probability and numerical algorithms

In this section, we discuss the mean exit time and the escape probability and their numerical schemes in order to quantify the dynamics of the stochastic differential equation (8).

## 4.1 Mean exit time

First, we recall the definition of the first exit time from the bounded domain $D = (x_1, x_3)$:

$$\tau(\omega) = \inf\{t > 0, X_t(\omega, x) \notin D\}.$$

The mean exit time $u(x) = \mathbb{E}\tau(\omega)$ then satisfies the following integro-differential equation [5] with an exterior boundary condition

$$\begin{aligned} Au(x) &= -1, \quad x \in D, \\ u &= 0, \quad x \in D^c, \end{aligned} \qquad (9)$$

where the generator $A$ of the solution process $X(t)$ is

$$\begin{aligned} Au &= f(x)u'(x) + \frac{d}{2}u''(x) \\ &+ \int_{\mathbb{R}^1 \setminus \{0\}} [u(x+y) - u(x) - \mathbf{1}_{\{|y|<1\}} yu'(x)] \, \nu_{\alpha,\beta}(\mathrm{d}y). \end{aligned} \qquad (10)$$

## 4.2 Escape probability

Now, we consider the escape probability of paths whose motion is described by Eq. (8). The likelihood that $X(t)$, starting at a point $x$, first exits from the domain $D$ by landing in the subset $E$ of $D^c$ is called the escape probability and is denoted as $p(x)$. This escape probability solves the following exterior Dirichlet problem [24]:

$$\begin{aligned} Ap(x) &= 0, \quad x \in D, \\ p &= \begin{cases} 1, & x \in E, \\ 0, & x \in D^c \setminus E, \end{cases} \end{aligned} \qquad (11)$$

where $A$ is the generator defined in (10).

We are concerned with how to enhance the likelihood of tumor extinction, so we choose $E = (-\infty, x_1]$, i.e., we examine the likelihood that the tumor state goes from $D = (x_1, x_3)$ (tumor) to $E = (-\infty, x_1]$ (tumor free).

## 4.3 Numerical algorithms

We only describe the numerical algorithms for Eq. (11), as the algorithm for Eq. (9) is similar. The algorithm below extends a numerical scheme in [11] for the case of symmetric Lévy noise to the case of asymmetric Lévy noise. For convenience, we use a general interval $D = (a, b)$, instead of $D = (x_1, x_3)$, in the following spatial discretization. Using (3), Eq. (11) can be rewritten as:

$$\frac{d}{2}p''(x) + f(x)p'(x) \tag{12}$$
$$+ \int_{\mathbb{R}^1\setminus\{0\}} [p(x+y) - p(x) - \mathbf{1}_{\{|y|<1\}}yp'(x)][\frac{C_1\mathbf{1}_{(0,+\infty)}(y)}{|y|^{1+\alpha}} + \frac{C_2\mathbf{1}_{(-\infty,0)}(y)}{|y|^{1+\alpha}}]dy = 0,$$

for $x \in (a,b)$; $p(x) = 1$ for $x \in (-\infty, a]$ and $p(x) = 0$ for $x \in [b, +\infty)$. Thus, we obtain:

$$\frac{d}{2}p''(x) + f(x)p'(x) + C_2 \int_{\mathbb{R}^1\setminus\{0\}} \frac{p(x+y) - p(x) - \mathbf{1}_{\{|y|<1\}}yp'(x)}{|y|^{1+\alpha}} dy \tag{13}$$
$$+ (C_1 - C_2) \int_{\mathbb{R}^{1+}} \frac{p(x+y) - p(x) - \mathbf{1}_{\{|y|<1\}}yp'(x)}{|y|^{1+\alpha}} dy = 0.$$

Because $\int_{\mathbb{R}^1\setminus\{0\}} \frac{\mathbf{1}_{\{|y|<1\}}yp'(x)}{|y|^{1+\alpha}} dy = \int_{\mathbb{R}^1\setminus\{0\}} \frac{\mathbf{1}_{\{|y|<\delta\}}yp'(x)}{|y|^{1+\alpha}} dy = 0$, we can replace the former by the latter in (13). To take care of the external condition, we divide the integral as $\int_{\mathbb{R}^1} = \int_{-\infty}^{a-x} + \int_{a-x}^{b-x} + \int_{b-x}^{+\infty}$ and choose $\delta = \min\{|a-x|, |b-x|\}$, and then we get the following formula:

$$\frac{d}{2}p''(x) + f(x)p'(x) + C_2 \int_{-\infty}^{a-x} \frac{p(x+y) - p(x) - \mathbf{1}_{\{|y|<\delta\}}yp'(x)}{|y|^{1+\alpha}} dy \tag{14}$$
$$+ C_2 \int_{b-x}^{\infty} \frac{p(x+y) - p(x) - \mathbf{1}_{\{|y|<\delta\}}yp'(x)}{|y|^{1+\alpha}} dy$$
$$+ C_2 \int_{a-x}^{b-x} \frac{p(x+y) - p(x) - \mathbf{1}_{\{|y|<\delta\}}yp'(x)}{|y|^{1+\alpha}} dy$$
$$+ (C_1 - C_2) \int_{\mathbb{R}^{1+}} \frac{p(x+y) - p(x) - \mathbf{1}_{\{|y|<1\}}yp'(x)}{|y|^{1+\alpha}} dy = 0.$$

By direct calculations, (14) can be further rewritten as:

$$\frac{d}{2}p''(x) + f(x)p'(x) - \frac{C_2}{\alpha}\left[\frac{1}{(x-a)^\alpha} + \frac{1}{(b-x)^\alpha}\right]p(x) \tag{15}$$
$$+ C_2 \int_{x-a}^{b-x} \frac{p(x+y) - p(x)}{|y|^{1+\alpha}} dy + C_2 \int_{a-x}^{x-a} \frac{p(x+y) - p(x) - yp'(x)}{|y|^{1+\alpha}} dy$$
$$+ (C_1 - C_2) \int_{\mathbb{R}^{1+}} \frac{p(x+y) - p(x) - \mathbf{1}_{\{|y|<1\}}yp'(x)}{|y|^{1+\alpha}} dy = -\frac{C_2}{\alpha}\left[\frac{1}{(x-a)^\alpha}\right],$$

for $x < \frac{a+b}{2}$, and

$$\frac{d}{2}p''(x) + f(x)p'(x) - \frac{C_2}{\alpha}\left[\frac{1}{(x-a)^\alpha} + \frac{1}{(b-x)^\alpha}\right]p(x) \qquad (16)$$
$$+C_2 \int_{a-x}^{x-b} \frac{p(x+y)-p(x)}{|y|^{1+\alpha}}\,dy + C_2 \int_{x-b}^{b-x} \frac{p(x+y)-p(x)-yp'(x)}{|y|^{1+\alpha}}\,dy$$
$$+(C_1-C_2)\int_{\mathbb{R}^{1+}} \frac{p(x+y)-p(x)-\mathbf{1}_{\{|y|<1\}}yp'(x)}{|y|^{1+\alpha}}\,dy = -\frac{C_2}{\alpha}\left[\frac{1}{(x-a)^\alpha}\right],$$

for $x \geq \frac{a+b}{2}$.

For the last term on the left hand of the Eq. (14), we have (Here, we omit the coefficient $(C_1 - C_2)$ and assume $b - a > 1$.):

$$\int_{\mathbb{R}^{1+}} \frac{p(x+y)-p(x)-\mathbf{1}_{\{|y|<1\}}yp'(x)}{|y|^{1+\alpha}}\,dy \qquad (17)$$
$$= \int_0^1 \frac{p(x+y)-p(x)-yp'(x)}{|y|^{1+\alpha}}\,dy + \int_1^{b-x}\frac{p(x+y)-p(x)}{|y|^{1+\alpha}}\,dy + \int_{b-x}^{+\infty}\frac{-p(x)}{|y|^{1+\alpha}}\,dy$$
$$= \int_0^1 \frac{p(x+y)-p(x)-yp'(x)}{|y|^{1+\alpha}}\,dy + \int_1^{b-x}\frac{p(x+y)-p(x)}{|y|^{1+\alpha}}\,dy - \frac{p(x)}{\alpha(b-x)^\alpha},$$

for $x \leq b-1$, i.e. $b-x \geq 1$, and

$$\int_{\mathbb{R}^+} \frac{p(x+y)-p(x)-\mathbf{1}_{\{|y|<1\}}yp'(x)}{|y|^{1+\alpha}}\,dy \qquad (18)$$
$$= \int_0^{b-x} \frac{p(x+y)-p(x)-yp'(x)}{|y|^{1+\alpha}}\,dy - \int_{b-x}^1 \frac{p(x)+yp'(x)}{|y|^{1+\alpha}}\,dy - \int_1^{+\infty}\frac{p(x)}{|y|^{1+\alpha}}\,dy$$
$$= \int_0^{b-x} \frac{p(x+y)-p(x)-yp'(x)}{|y|^{1+\alpha}}\,dy - \int_{b-x}^1 \frac{p(x)+yp'(x)}{|y|^{1+\alpha}}\,dy - \frac{p(x)}{\alpha},$$

for $x > b-1$, i.e. $b-x < 1$.

Now, let us take the appropriate stepsize $h$, so that $\frac{1}{h}, \frac{a}{h}, \frac{b}{h}, \frac{a+b}{2h}$ are integers, and define $x_j = jh$ for $\frac{a-b}{h} \leq j \leq \frac{b-a}{h}$. We use the notation $P_j$ to indicate the numerical solution of $p$ at $x_j$. Then, the first four terms on the left hand of the Eqs. (15) and (16) can be discretized, respectively, by the central difference scheme for derivatives and "punched-hole" trapezoidal rule:

$$\frac{d}{2}\frac{P_{j-1}-2P_j+P_{j+1}}{h^2} + f(x_j)\frac{P_{j+1}-P_{j-1}}{2h} - \frac{C_2}{\alpha}\left[\frac{1}{(x_j-a)^\alpha} + \frac{1}{(b-x_j)^\alpha}\right]P_j$$
$$+C_2 h \sum_{k=j-\frac{a}{h}}^{\frac{b}{h}-j}{}'' \frac{P_{j+k}-P_j}{|x_k|^{1+\alpha}} + C_2 h \sum_{k=\frac{a}{h}-j, k\neq 0}^{j-\frac{a}{h}}{}'' \frac{P_{j+k}-P_j-(P_{j+1}-P_{j-1})x_k/2h}{|x_k|^{1+\alpha}},$$
$$\qquad (19)$$

where $j = \frac{a}{h}+1, \frac{a}{h}+2, \cdots, \frac{a+b}{2h}-1$. The meaning of the modified summation symbol $\sum''$ is that the quantities corresponding to the two endpoints of the integral interval should be multiplied by $1/2$.

$$\frac{d}{2}\frac{P_{j-1}-2P_j+P_{j+1}}{h^2} + f(x_j)\frac{P_{j+1}-P_{j-1}}{2h} - \frac{C_2}{\alpha}\left[\frac{1}{(x_j-a)^\alpha}+\frac{1}{(b-x_j)^\alpha}\right]P_j$$
$$+ C_2 h \sum_{k=\frac{a}{h}-j}^{j-\frac{b}{h}}{}'' \frac{P_{j+k}-P_j}{|x_k|^{1+\alpha}} + C_2 h \sum_{k=j-\frac{b}{h},k\neq 0}^{\frac{b}{h}-j}{}'' \frac{P_{j+k}-P_j-(P_{j+1}-P_{j-1})x_k/2h}{|x_k|^{1+\alpha}},$$
(20)

where $j = \frac{a+b}{2h}, \frac{a+b}{2h}+1, \cdots, \frac{b}{h}-1$.

For the last term on the left hand of the Eqs. (15) and (16), we use the upwind finite difference scheme. That is, when $-(C_1-C_2)\int_{\mathbb{R}^{1+}}\frac{\mathbf{1}_{\{|y|<1\}}y}{|y|^{1+\alpha}}\,dy > 0$, we use the forward finite difference scheme to discretize $p'(x)$, while when $-(C_1-C_2)\int_{\mathbb{R}^{1+}}\frac{\mathbf{1}_{\{|y|<1\}}y}{|y|^{1+\alpha}}\,dy \leq 0$, we use the backward finite difference scheme to discretize $p'(x)$, i.e., as in [7],

$$p'(x) = \begin{cases} \dfrac{p(x+h)-p(x)}{h}, & -(C_1-C_2)\int_{\mathbb{R}^{1+}}\dfrac{\mathbf{1}_{\{|y|<1\}}y}{|y|^{1+\alpha}}\,dy > 0, \\ \dfrac{p(x)-p(x-h)}{h}, & -(C_1-C_2)\int_{\mathbb{R}^{1+}}\dfrac{\mathbf{1}_{\{|y|<1\}}y}{|y|^{1+\alpha}}\,dy \leq 0. \end{cases} \quad (21)$$

So the formula (17) can be discretized as:

$$-\frac{1}{\alpha(b-x_j)^\alpha}P_j + h\sum_{k=\frac{1}{h}}^{\frac{b}{h}-j}{}''\frac{P_{j+k}-P_j}{|x_k|^{1+\alpha}} + h\sum_{k=0,k\neq 0}^{\frac{1}{h}}{}''\frac{P_{j+k}-P_j-(P_{j+1}-P_j)x_k/h}{|x_k|^{1+\alpha}},$$
(22)

for $\beta < 0$, and

$$-\frac{1}{\alpha(b-x_j)^\alpha}P_j + h\sum_{k=\frac{1}{h}}^{\frac{b}{h}-j}{}''\frac{P_{j+k}-P_j}{|x_k|^{1+\alpha}} + h\sum_{k=0,k\neq 0}^{\frac{1}{h}}{}''\frac{P_{j+k}-P_j-(P_j-P_{j-1})x_k/h}{|x_k|^{1+\alpha}},$$
(23)

for $\beta \geq 0$. By the same reason, we can get the discretization for the formula (18):

$$-\frac{1}{\alpha}P_j + h\sum_{\frac{b}{h}-j}^{\frac{1}{h}}{}''\frac{-P_j-(P_{j+1}-P_j)x_k/h}{|x_k|^{1+\alpha}} + h\sum_{k=0,k\neq 0}^{\frac{b}{h}-j}{}''\frac{P_{j+k}-P_j-(P_{j+1}-P_j)x_k/h}{|x_k|^{1+\alpha}},$$
(24)

for $\beta < 0$, and

$$-\frac{1}{\alpha}P_j + h\sum_{\frac{b}{h}-j}^{\frac{1}{h}}{''}\frac{-P_j - (P_j - P_{j-1})x_k/h}{|x_k|^{1+\alpha}} + h\sum_{k=0,k\neq 0}^{\frac{b}{h}-j}{''}\frac{P_{j+k} - P_j - (P_j - P_{j-1})x_k/h}{|x_k|^{1+\alpha}}, \tag{25}$$

for $\beta \geq 0$. The boundary conditions require that $P_j = 1$ for $j \leq \frac{a}{h}$ and $P_j = 0$ for $j \geq \frac{b}{h}$.

The right hand of the Eq. (15) can be discreted as:

$$-\frac{C_2}{\alpha}[\frac{1}{(x_j - a)^\alpha}]. \tag{26}$$

Thus, we can obtain the escape probability $p(x)$, $x \in D$, by solving the equations (19)-(26).

A numerical scheme to solve Eq. (9) for mean exit time $u(x)$ is similar.

# 5 Numerical results

We fix the chemical system parameters $\theta = 0.1, \gamma = 3.0$ as suggested in [10], and focus on the impact of asymmetric Lévy noise on the mean exit time $u$ and the escape probability $p$, in order to gain understanding of the tumor evolution under uncertainty. With these system parameters, the stable states of the deterministic tumor growth model are $x_1 = 0$ and $x_3 = 5$. The interval $D = (0, 5)$ encloses the tumor cell density from 0 (tumor-free state) to 5 (stable tumor state).

## 5.1 Mean exit time

We compute the mean exit time $u(x)$ of the tumor density $X_t$ from $x \in D = (0, 5)$, between the tumor-free state and the stable tumor state. For example, $u(3)$ is the mean time that the tumor cell density, starting at the initial density 3, remains within the range $D$, before 'exiting' to outside $D$. It quantifies how long the tumor cell density stays between 0 (tumor-free state) and 5 (stable tumor state).

In Figures 3 and 4, we observe that when the stability index $\alpha$ is fixed approximately below 1.77, for the inchoate patients, smaller skewness parameter $\beta$ leads to shorter mean exit time. While for the advanced patients, the situation becomes opposite. This suggests that asymmetry in the noise (described by $\beta$) plays a crucial role in the time span a patient remain in the tumor state, and this role is opposite for inchoate and advanced patients.

Besides, when the stability index $\alpha \geq 1.77$, for all the patients, the mean exit time $u$ decreases as the skewness parameter $\beta$ decreases. From the Figures 5 and 6, we find that when the skewness parameter $\beta$ is fixed, the behavior of the mean exit time under various stability indices $\alpha$ is similar to the symmetric case (Figures 5(c) and 6(c) for $\beta = 0$). That is, at the early ($x$ near 0) and advanced ($x$ close to 5) tumor stages, the mean exit time decreases as the stability index $\alpha$ increases, while in the middle tumor stage, this relationship becomes reversed. Moreover, comparing Figure 3 with Figure 4, and also Figure 5 with Figure 6, we notice that the input of the Gaussian noise component (i.e., the diffusion coefficient $d > 0$) just shortens the mean exit time, but appears not to change the evolution time essentially.

Figure 11 is a 3D plot of the above numerical results.

## 5.2 Escape probability

We like to do further investigation on the tumor evolution after the system exits the interval $D = (0, 5)$. Will it reach the tumor-free state (exit to the left of $D = (0, 5)$) ? We compute the escape probability $p(x)$, i.e., the likelihood that the tumor cell density $x \in (0, 5)$ to become tumor-free ($x = 0$). See Figures 7 – 10. In these figures, high escape probability $p$ values indicate high likelihood for a tumor to become tumor-free.

In Figures 7 and 8, we observe that when the stability index $\alpha$ is fixed, either with or without Gaussian noise component (i.e., diffusion coefficient $d = 0$ or $d = 1$), the escape probability $p$ increases as the skewness parameter $\beta$ decreases. Besides, the smaller the stability index $\alpha$ is, the more effects the skewness parameter $\beta$ has on the escape probability. As the stability index $\alpha$ becomes larger (close to 2), for all the skewness parameter $\beta$, the escape probability tends to the symmetric case (i.e. $\beta = 0$). Especially, when $\alpha = 1.9$ (see Figures 7(e) and 8(e)), the escape probability does not change much with the skewness parameter $\beta$.

As seen in Figures 9 and 10, when the skewness parameter $\beta$ is near the extremes, the escape probability increases as the stability index $\alpha$ either decreases (in the $\beta \leq -0.99$ case) or increases (in the $\beta \geq 0.99$ case). However, when $-0.99 < \beta < 0.99$, due to the competition between the stability index $\alpha$ and the skewness parameter $\beta$, the behavior of the escape probability does not change monotonically with $\alpha$. Instead, they are similar to the symmetric case (see Figures 9(c) and 10(c) for $\beta = 0$ case), in which case there exists a critical point. When the density of the tumor cells is less than the critical point, the escape probability increases with the increasing $\alpha$. While when the density of the tumor cells is more than the critical point, the escape probability shows an opposite trend with the increasing $\alpha$.

Comparing all these figures, we find that the smaller the stability index $\alpha$

and the skewness parameter $\beta$ are, the shorter the mean exit time $u$ becomes and the higher the escape probability $p$ grows. Especially, when $\alpha = 0.1$ and $\beta = -1.0$, no matter at what tumor stage (either the early stage or the advanced stage) the patients are, the probability of curing the tumor is almost surely ($p \sim 1$, see Figure 7(a) for $\beta = -1$). This is different from the symmetric noise case, which suggests that we should make the therapy strategy according to the stage of the patients (see Figures 9(c) and 10(c)).

Figure 12 is a 3D plot of the above numerical results.

# 6  Conclusions

In this paper, we have examined the impact of asymmetric non-Gaussian environmental fluctuations on the dynamical evolution of a tumor growth model with immunization. The mean exit time and the escape probability are computed to quantify the mean lifetime the tumor cells remains between the tumor-free state and the stable tumor state, and the likelihood that tumor with a certain initial density becomes extinct (i.e., tumor-free). These two quantities are described by exterior boundary value problems involving nonlocal operators. By the numerical experiments, we find that the parameters of asymmetric non-Gaussian Lévy noise have significant influences on the mean lifetime and the escape probability. Especially, the skewness parameter plays an important role in controlling the tumor cells overall evolution. By choosing the appropriate skewness parameter, we observe that it is likely to slow the tumor progression and at the same time, enhance the likelihood to induce tumor extinction.

# Acknowledgement

We would like to thank Xiaofan Li and Ting Gao for helpful discussions on the numerical scheme.

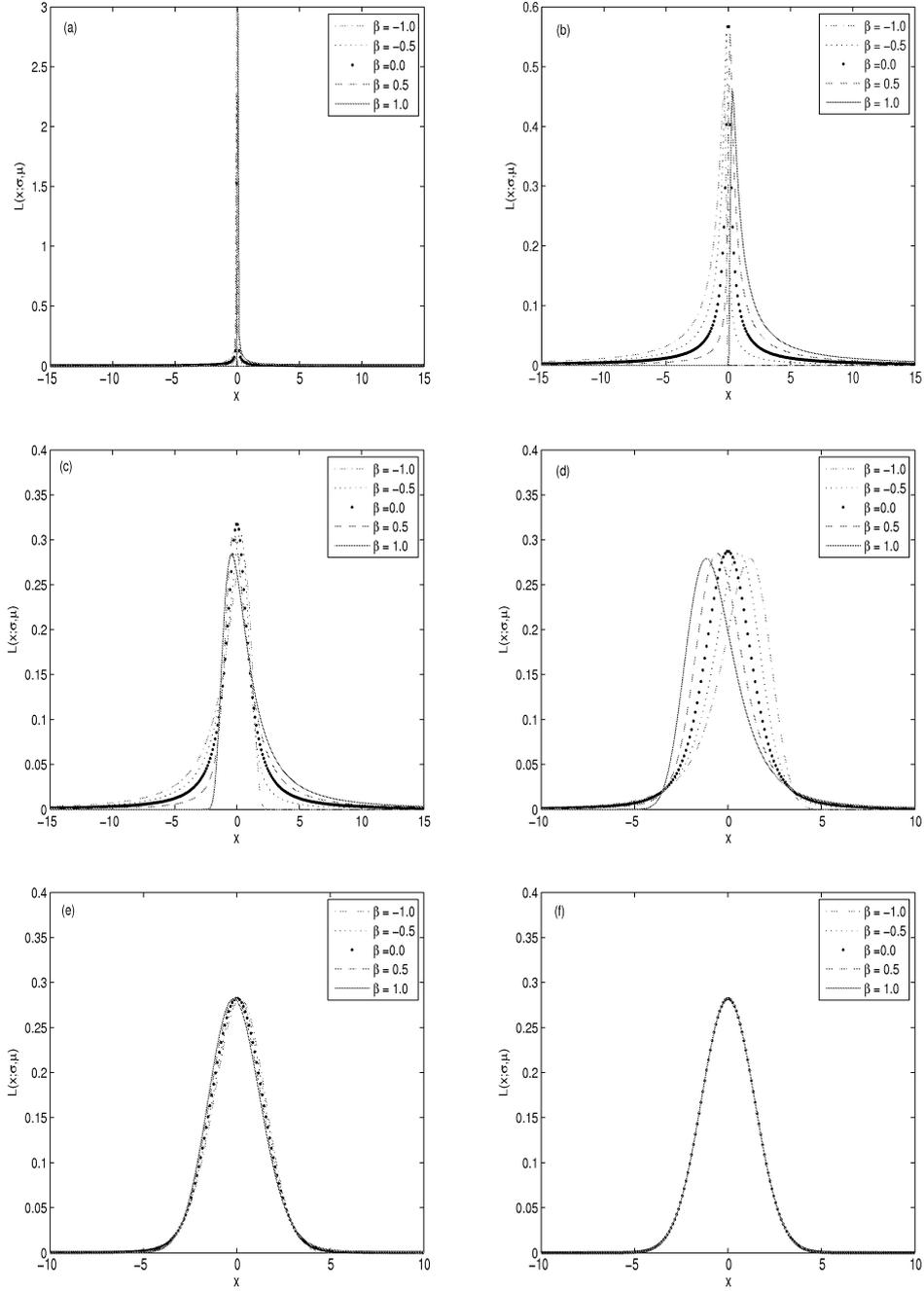

Figure 1: Probability density functions for $L_1 \sim S_\alpha(1, \beta, 0)$ for different values of $\alpha$ and $\beta$: (a) $\alpha = 0.1$; (b) $\alpha = 0.5$; (c) $\alpha = 1$; (d) $\alpha = 1.5$; (e) $\alpha = 1.9$; (f) $\alpha = 2$.

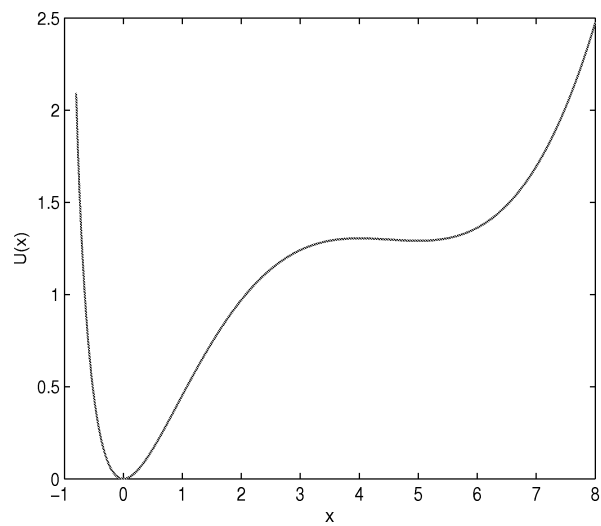

Figure 2: The potential function $U(x)$ for $\theta = 0.1, \gamma = 3.0$. The two minima represent the tumor-free state at $x_1 = 0$ and the stable tumor state at $x_3 = 5$, respectively.

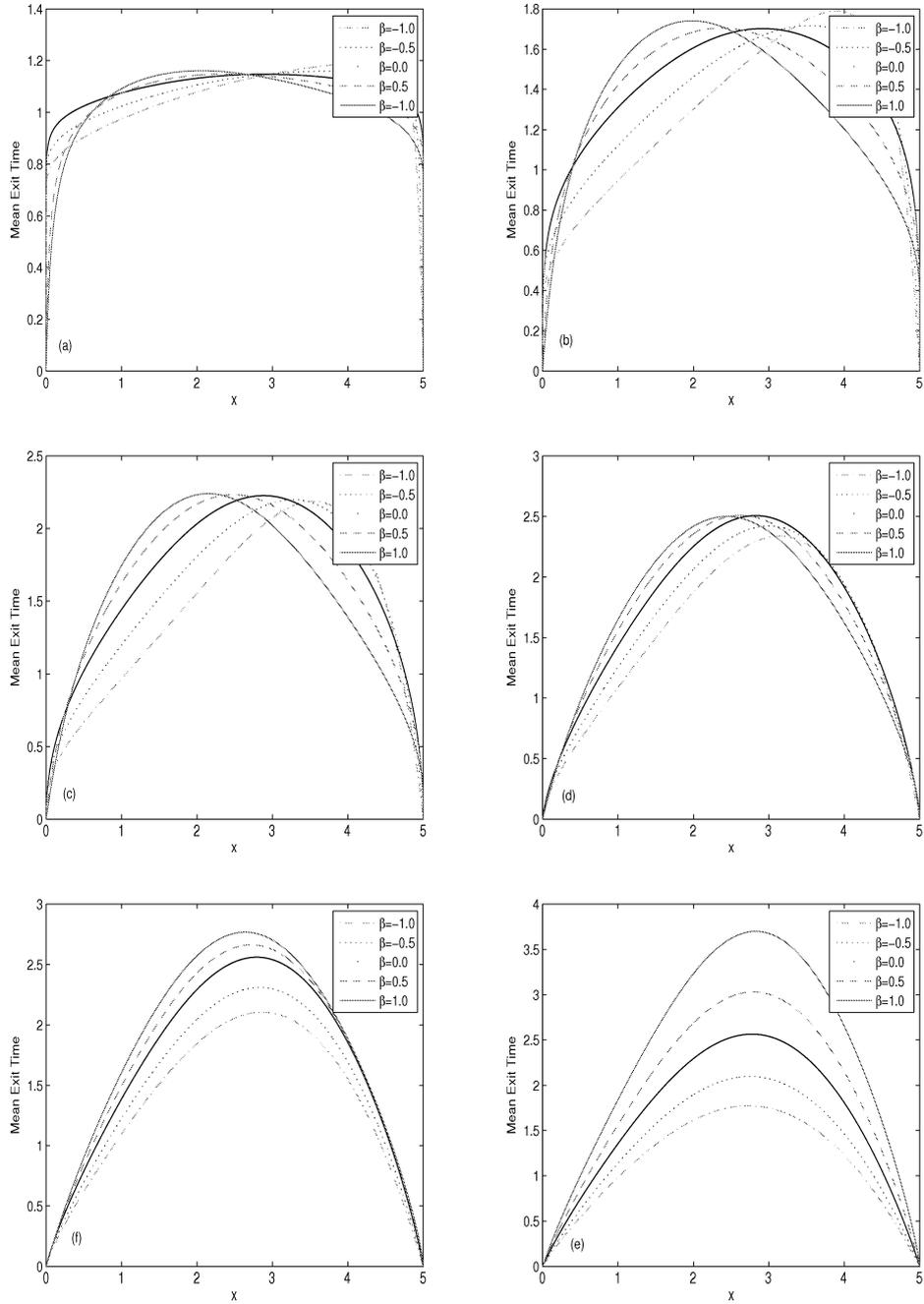

Figure 3: Mean exit time $u(x)$ with pure $\alpha$-satble Lévy noise (i.e. $d = 0$): (a) $\alpha = 0.1$; (b) $\alpha = 0.5$; (c) $\alpha = 1.0$; (d) $a = 1.5$; (e) $\alpha = 1.77$; (f) $\alpha = 1.9$.

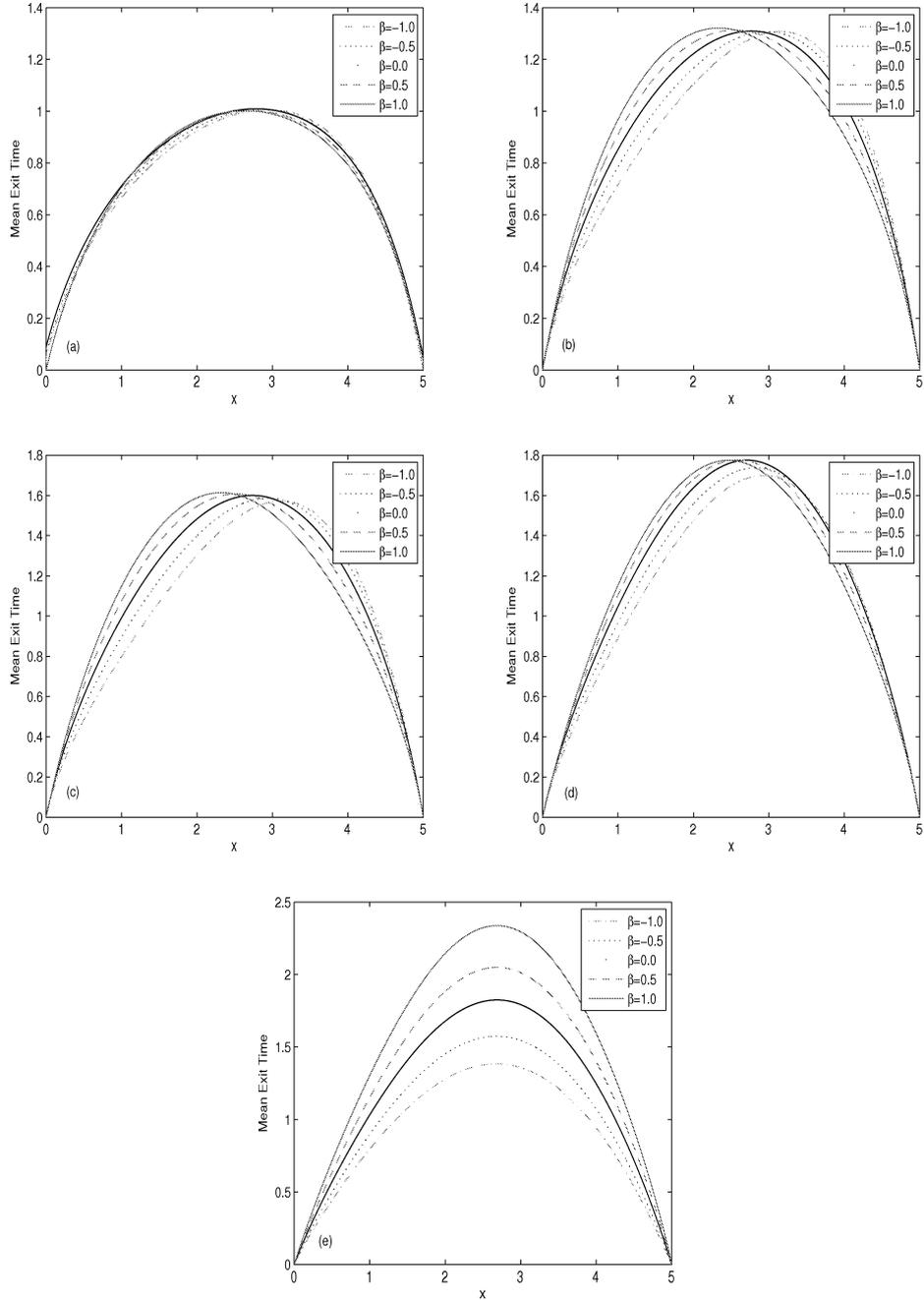

Figure 4: Mean exit time $u(x)$ when pure $\alpha$-satble Lévy noise is combined with Gaussian noise (i.e. $d = 1.0$): (a) $\alpha = 0.1$; (b) $\alpha = 0.5$; (c) $\alpha = 1.0$; (d) $\alpha = 1.5$; (e) $\alpha = 1.9$.

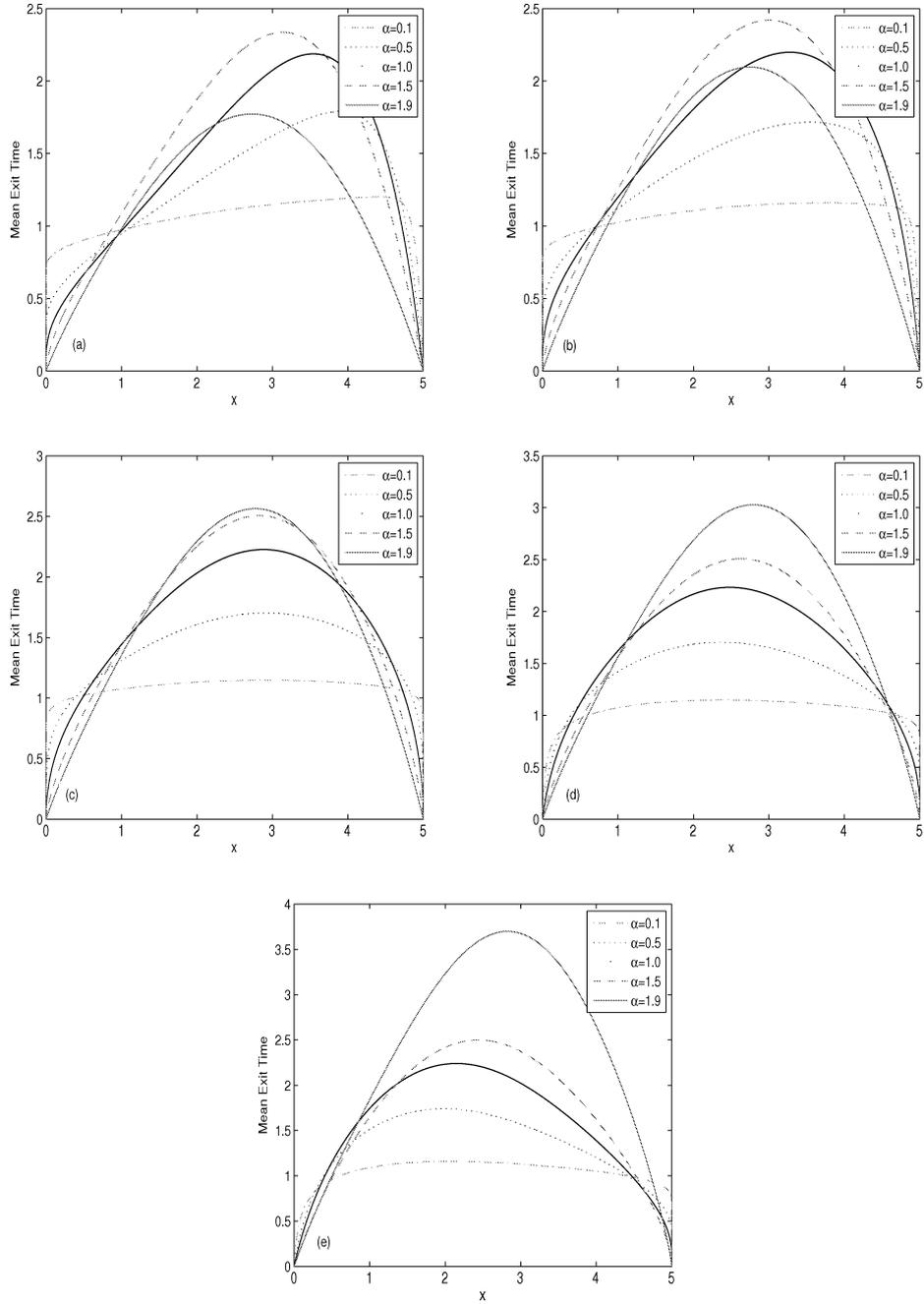

Figure 5: Mean exit time $u(x)$ with pure $\alpha$-satble Lévy noise (i.e. $d = 0$): (a) $\beta = -1.0$; (b) $\beta = -0.5$; (c) $\beta = 0.0$; (d) $\beta = 0.5$; (e) $\beta = 1.0$.

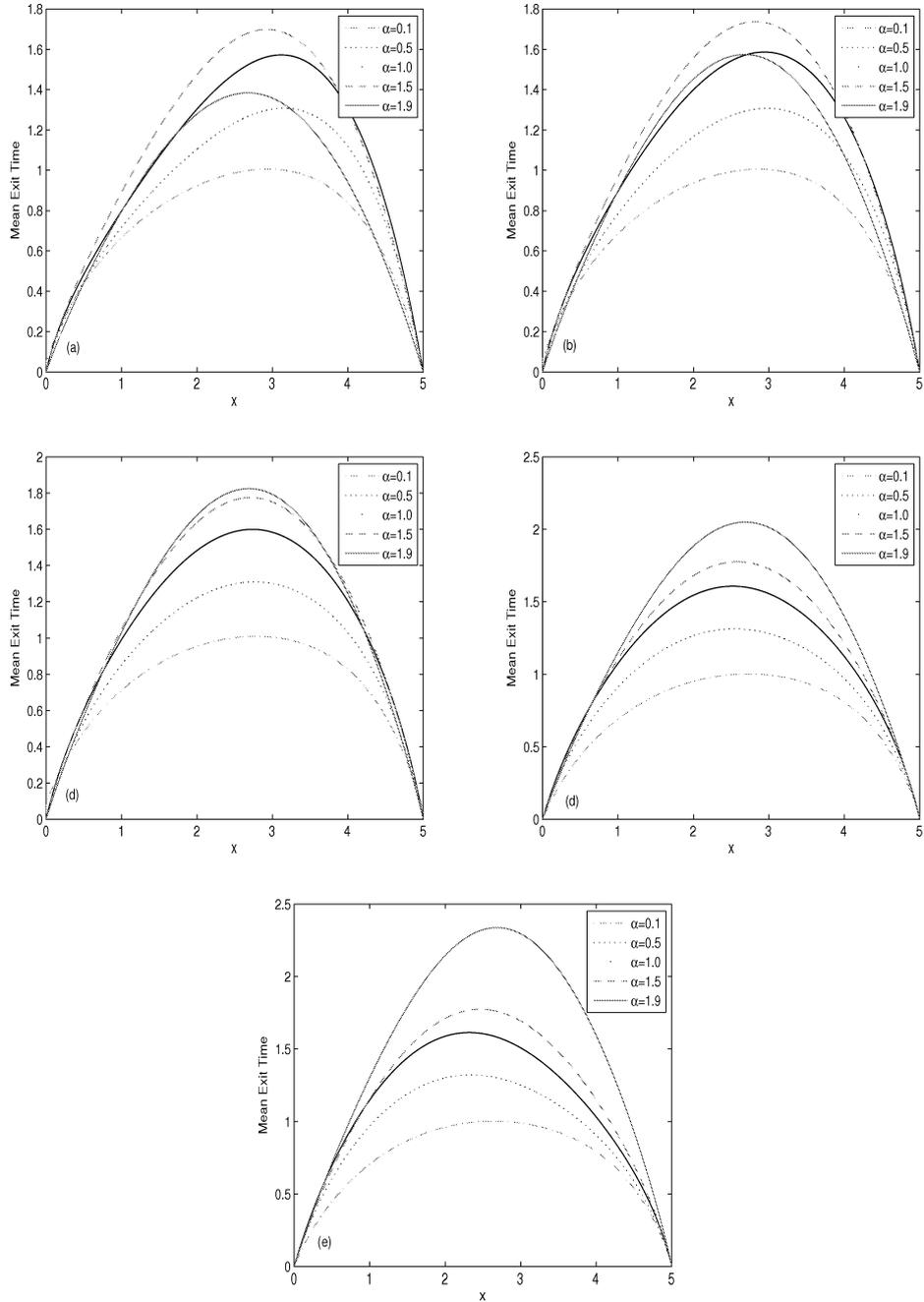

Figure 6: Mean exit time $u(x)$ when pure $\alpha$-satble Lévy noise is combined with Gaussian noise (i.e. $d = 1.0$): (a) $\beta = -1.0$; (b) $\beta = -0.5$; (c) $\beta = 0.0$; (d) $\beta = 0.5$; (e) $\beta = 1.0$.

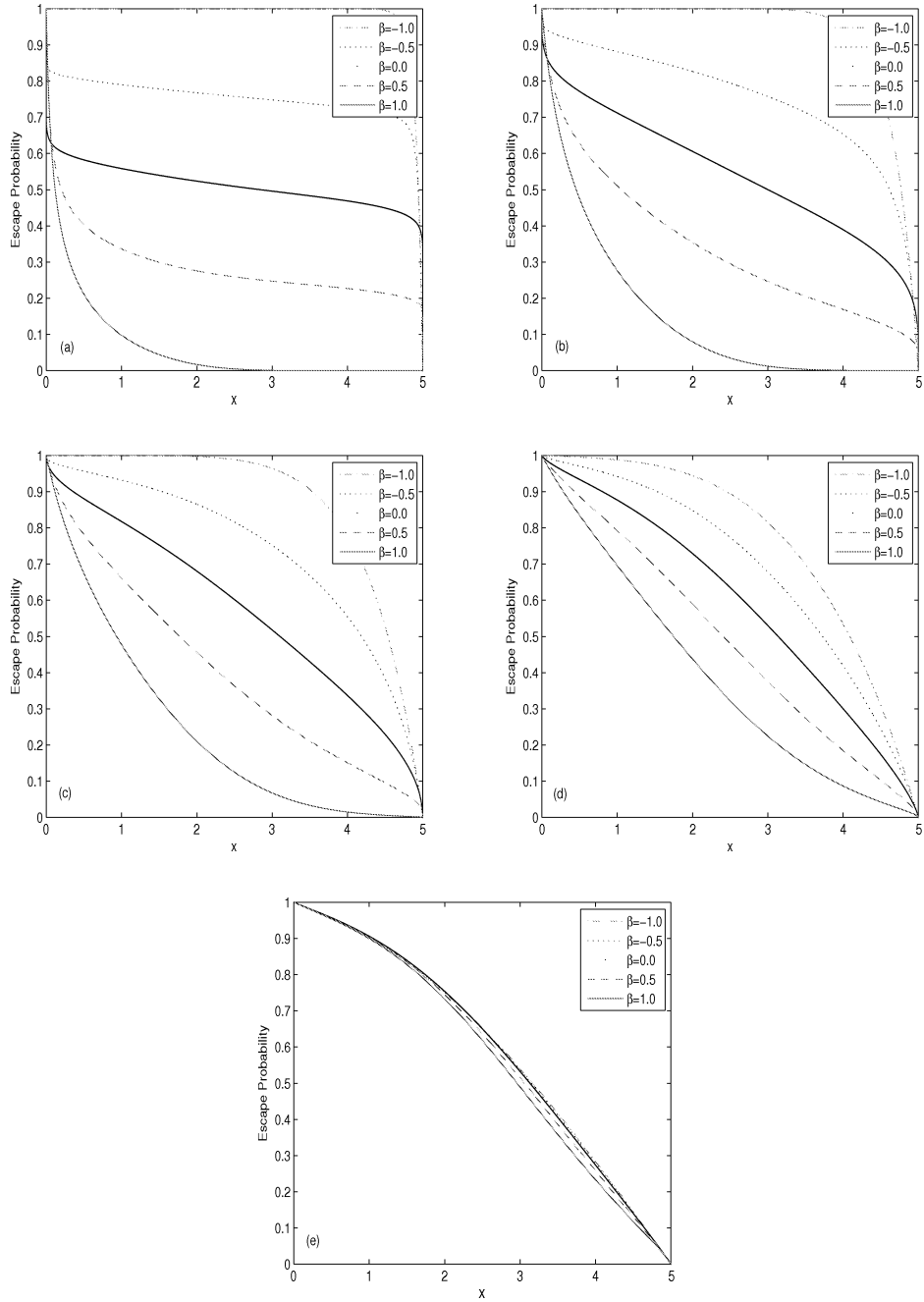

Figure 7: Escape Probability $p(x)$ with pure $\alpha$-satble Lévy noise (i.e. $d = 0$): (a) $\alpha = 0.1$; (b) $\alpha = 0.5$; (c) $\alpha = 1.0$; (d) $\alpha = 1.5$; (e) $\alpha = 1.9$.

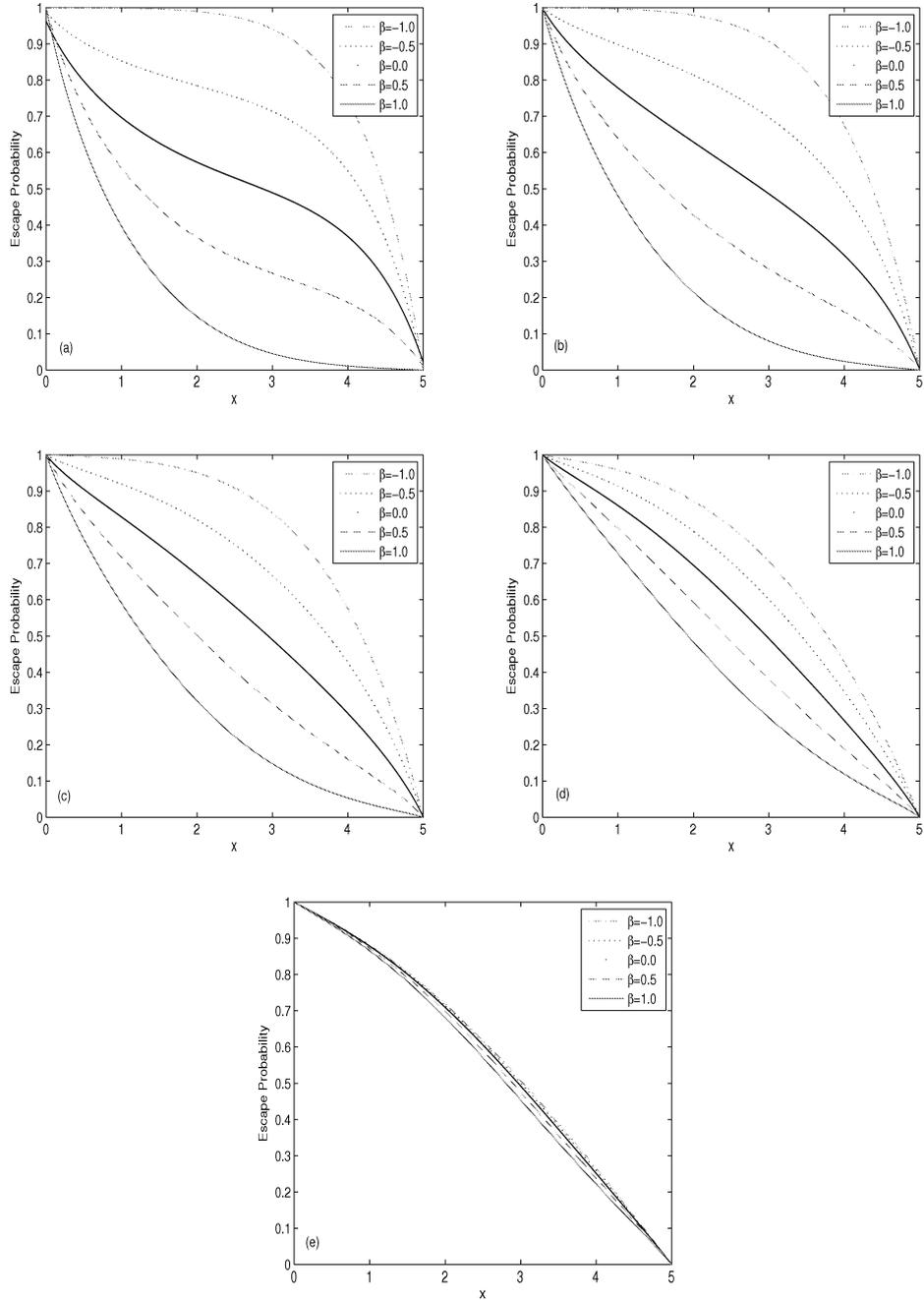

Figure 8: Escape Probability $p(x)$ when pure $\alpha$-satble Lévy noise is combined with Gaussian noise (i.e. $d = 1.0$): (a) $\alpha = 0.1$; (b) $\alpha = 0.5$; (c) $\alpha = 1.0$; (d) $a = 1.5$; (e) $\alpha = 1.9$.

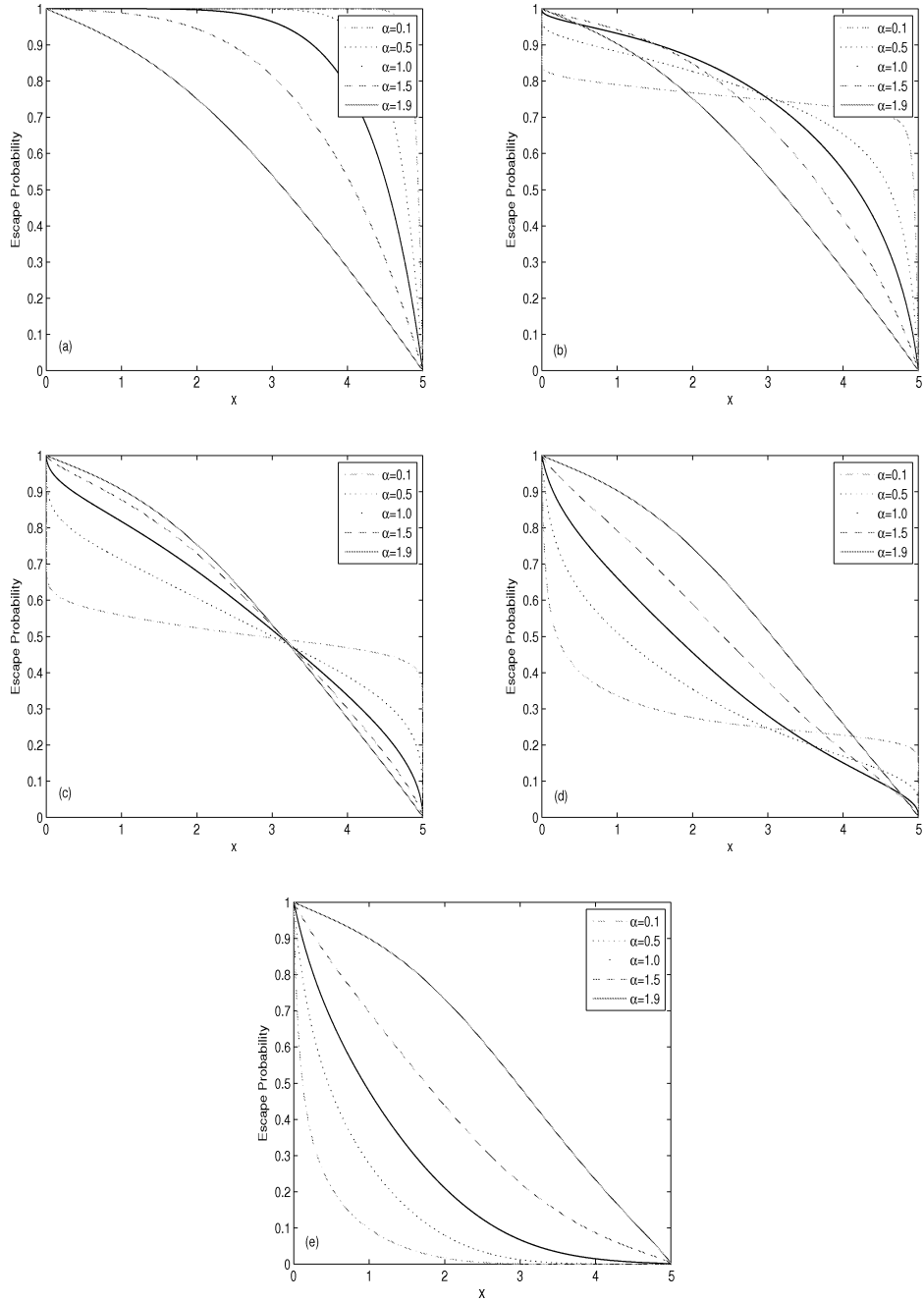

Figure 9: Escape Probability $p(x)$ with pure $\alpha$-satble Lévy noise (i.e. $d = 0$): (a) $\beta = -1.0$; (b) $\beta = -0.5$; (c) $\beta = 0.0$; (d) $\beta = 0.5$; (e) $\beta = 1.0$.

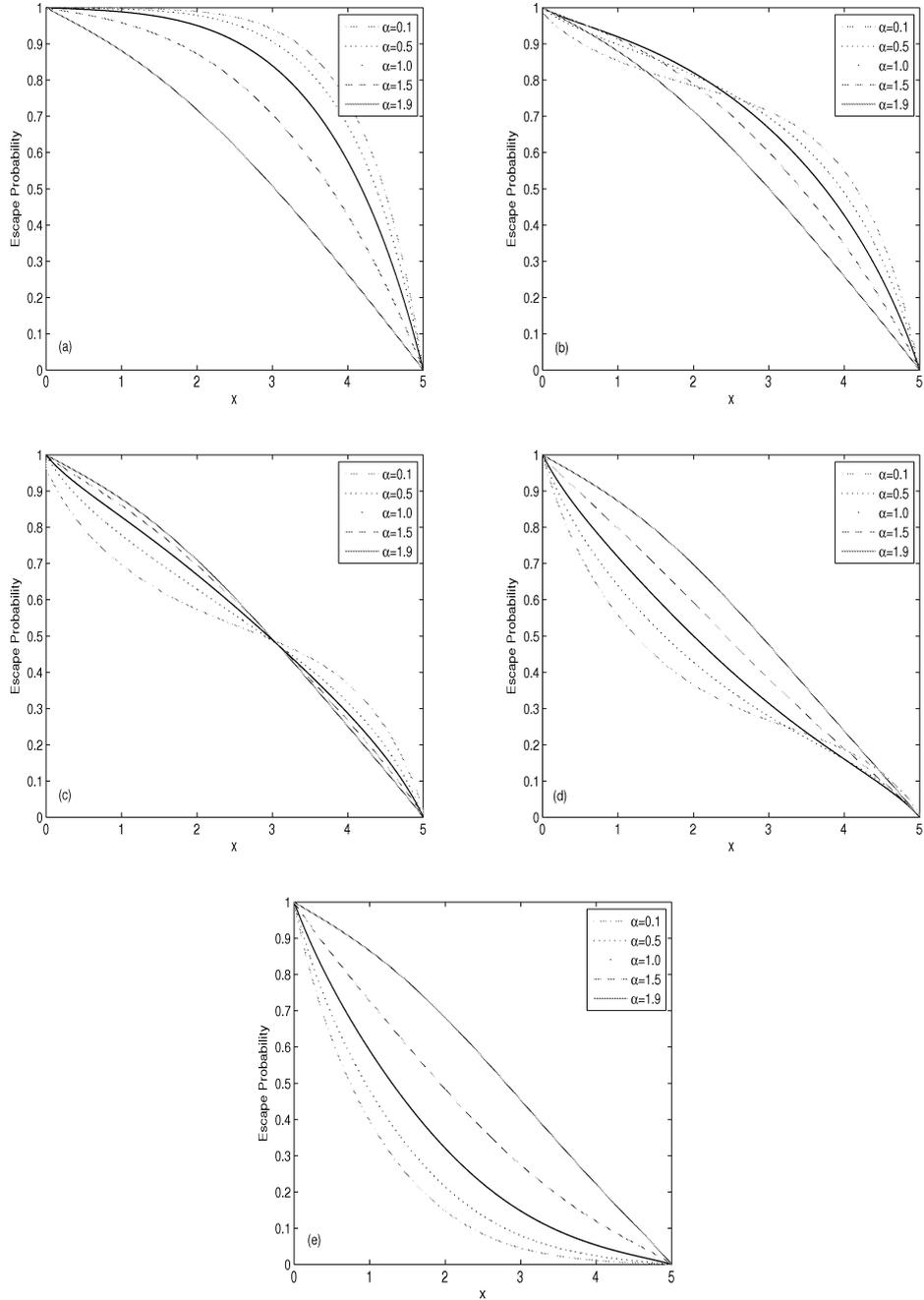

Figure 10: Escape Probability $p(x)$ when pure $\alpha$-satble Lévy noise is combined with Gaussian noise (i.e. $d = 1.0$): (a) $\beta = -1.0$; (b) $\beta = -0.5$; (c) $\beta = 0.0$; (d) $\beta = 0.5$; (e) $\beta = 1.0$.

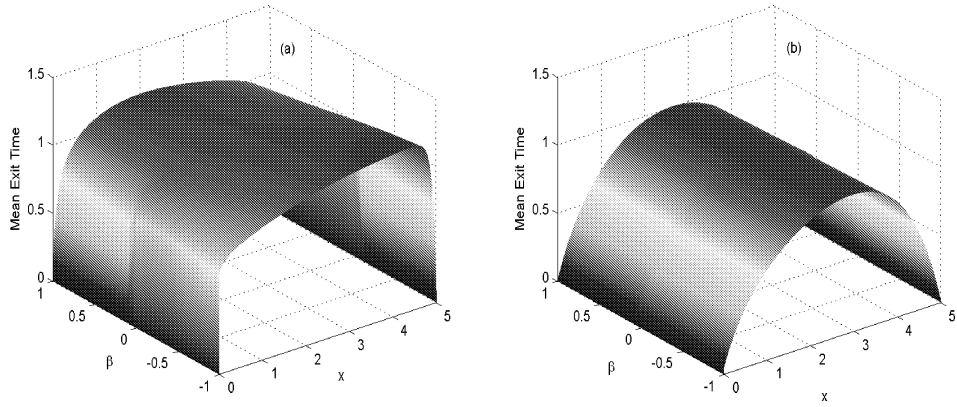

Figure 11: Mean exit time $u(x)$ in 3-dimension plane: (a) $d = 0, \alpha = 0.1$; (b) $d = 1, \alpha = 0.1$;

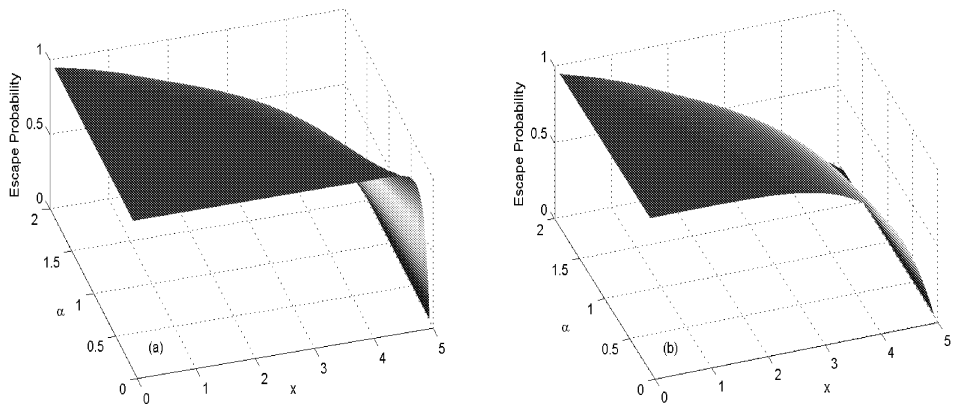

Figure 12: Escape Probability $p(x)$ in 3-dimension plane: (a) $d = 0, \beta = -1.0$; (b) $d = 1.0, \beta = -1.0$;